\documentclass[sn-mathphys-num]{sn-jnl}% Math and Physical Sciences Numbered Reference Style 
\usepackage{lmodern}
\usepackage{graphicx}%
\usepackage{multirow}%
\usepackage{amsmath,amssymb,amsfonts}%
\usepackage{amsthm}%
\usepackage[title]{appendix}%
\usepackage{xcolor}%
\usepackage{textcomp}%
\usepackage{manyfoot}%
\usepackage{booktabs}%
\usepackage{algorithm}%
\usepackage{algorithmicx}%
\usepackage{algpseudocode}%
\usepackage{listings}%
\newtheorem{theorem}{Theorem}%  meant for continuous numbers
\begin{document}

\title[Evaluation codes arising from symmetric polynomials]{Evaluation codes arising from symmetric polynomials}

%%=============================================================%%
%% GivenName	-> \fnm{Joergen W.}
%% Particle	-> \spfx{van der} -> surname prefix
%% FamilyName	-> \sur{Ploeg}
%% Suffix	-> \sfx{IV}
%% \author*[1,2]{\fnm{Joergen W.} \spfx{van der} \sur{Ploeg} 
%%  \sfx{IV}}\email{iauthor@gmail.com}
%%=============================================================%%

\author[1,2]{\fnm{Barbara} \sur{Gatti}}\email{barbara.gatti@unisalento.it}
\equalcont{These authors contributed equally to this work.}

\author[2]{\fnm{G\'abor} \sur{Korchm\'aros}}\email{gabor.korchmaros@unibas.it}
\equalcont{These authors contributed equally to this work.}

\author[3]{\fnm{G\'abor P.} \sur{Nagy}}\email{nagyg@math.u-szeged.hu}
\equalcont{These authors contributed equally to this work.}

\author[4]{\fnm{Vincenzo} \sur{Pallozzi Lavorante}}\email{vincenzop@usf.edu}
\equalcont{These authors contributed equally to this work.}

\author[1,2]{\fnm{Gioia} \sur{Schulte}}\email{gioia.schulte@unisalento.it}
\equalcont{These authors contributed equally to this work.}

\affil[1]{\orgdiv{Department of Mathematics and Physics ”Ennio de Giorgi”}, \orgname{University of Salento}, \orgaddress{\street{Via per Arnesano
}, \city{Lecce}, \postcode{73100 }, \state{Italy}}}

\affil[2]{\orgdiv{Department of Mathematics, Computer Science and Economics}, \orgname{University of Basilicata}, \orgaddress{\street{Contrada Macchia Romana}, \city{Potenza}, \postcode{85100}, \state{Italy}}}

\affil[3]{\orgdiv{Bolyai Institute
}, \orgname{University of Szeged
}, \orgaddress{\street{Aradi vértan\'uk tere 1
}, \city{Szeged}, \postcode{H-6720}, \state{Hungary}}}

\affil[4]{\orgdiv{Department of Mathematics and Statistics}, \orgname{University of South Florida}, \city{Tampa}, \postcode{33620}, \state{Florida}, \country{USA}}

%%==================================%%
%% Sample for unstructured abstract %%
%%==================================%%

\abstract{Datta and Johnsen (Des. Codes and Cryptogr., {\bf{91}} (2023), 747-761) introduced a new family of evalutation codes in an affine space of dimension $\ge 2$ over a finite field $\mathbb{F}_q$ where linear combinations of elementary symmetric polynomials are evaluated on the set of  all points with pairwise distinct coordinates. In this paper, we propose a generalization by taking low dimensional linear systems of symmetric polynomials. Computation for small values of $q=7,9$ shows that carefully chosen generalized Datta-Johnsen codes $\left[\frac{1}{2}q(q-1),3,d\right]$ have minimum distance $d$ equal to the optimal value minus 1.}


%%\pacs[JEL Classification]{D8, H51}


\maketitle

\section{Introduction}\label{sec1}

Evaluation codes are linear codes which are obtained by evaluating polynomials in $\mathbb{F}_q[X_1,\ldots,X_m]$ on some set of points of the $m$-dimensional affine space $AG(m,q)$ defined over the finite field $\mathbb{F}_q$ of order $q$.  The various known constructions of evaluation codes differ in how the polynomials and points are chosen, and include Reed-Solomon codes, Reed-Muller codes, monomial codes, Cartesian codes, and toric codes. The unique constraint on the choice of such polynomials is that they must form a finite dimensional $\mathbb{F}_q$-subspace of $\mathbb{F}_q[X_1,\ldots,X_m]$. This occurs, for instance, when all symmetric polynomials of bounded degrees in $\mathbb{F}_q[X_1,\ldots,X_m]$ are taken. Very recently, Datta and Johnsen  \cite{datta2023codes} considered the case of all $\mathbb{F}_q$-linear combinations of elementary symmetric polynomials in $\mathbb{F}_q[X_1,\ldots,X_m]$, and investigated the associated evaluation code on the set of all distinguished points in $AG(m,q)$, i.e. on all points with pairwise distinct coordinates. They computed its relevant parameters which are $$n=\binom{q}{m}m!,\,k=m+1, d=(q-m)\binom{q-1}{m-1}(m-1)!,$$ and observed that its relative minimum distance is the same as that of Reed-Muller codes. However, they found that the relative dimension (or rate) of their code is not as good. For this reason, Datta and Jonhsen exhibited a modified version, here named reduced Datta-Johnsen code, that has the same relative minimum distance, but a better rate. The proposed modification consists of evaluating symmetric polynomials on an ordered set $\mathcal{Q}$ of representatives of equivalence classes of distinguished points of $AG(m,q)$  where two distinct distinguished points are equivalent when they have the same coordinates but in a different order. For $m < q$, the reduced Datta-Johnsen code $C_m'$ is a non-degenerate $[N, K, D]$ linear code with $$N = \binom{q}{m},\, K=m+1,\, D=\binom{q}{m}- \binom{q-1}{m-1}.$$
%where the symmetric polynomials on the set of unramified points of the quotient variety $\mathbb{F}_q^m/{\rm{Sym}}_m$; see section \ref{av}.
The Datta-Johnsen codes are sub-codes of Reed-Muller codes and generated by minimum weight codewords.

A generalization of the Datta-Johnsen codes is found in \cite{Mi2} where the approach is a combination of Galois theoretical methods with Weil-type bounds for hypersurfaces.

Another generalization of the Datta-Johnsen code arises from any finite dimensional $\mathbb{F}_q$-subspace $V$ consisting of symmetric polynomials in $\mathbb{F}_q[X_1,\ldots,X_m]$, by evaluating those symmetric polynomials on the set of all distinguished points in $AG(m,q)$. In this context, the {\emph{reduced generalized Datta-Johnsen code}} is the evaluation code obtained by evaluating polynomials in $V$ on the set $\mathcal{Q}$.
The Datta-Johnsen code and its reduction correspond to the simplest case, i.e. $\dim(V)=1$. In this paper, we work out the case for $m=2$ and for certain subspaces $V$ of $\dim(V)\le 6$.

To work with linear systems of symmetric polynomials, our essential tool is the map $\Phi_m$ which takes a polynomial $f\in\mathbb{F}_q[X_1,\ldots,X_m]$ to the polynomial $\Phi_m(f)=f(\sigma_m^1(x),\dots,\sigma_m^m(x))$  where $\sigma_m^i(x)$ denotes the $i$-th elementary symmetric polynomial. From the fundamental theorem of symmetric polynomials,
${\rm{Im}}(\Phi_m)$ consists of all symmetric polynomials in $\mathbb{F}_q[X_1,\ldots,X_m]$.
Since $\Phi_m$ is an $\mathbb{F}_q$-linear map, every $\mathbb{F}_q$-subspace of symmetric polynomials is the image of a (unique) $\mathbb{F}_q$-subspace of $\mathbb{F}_q[X_1,\ldots,X_m]$. In particular, the linear system of all polynomials $\Sigma_m(t)$ of degree $\le t$, as well as any linear subsystem of $\Sigma_m(t)$, gives rise to an $\mathbb{F}_q$-subspace of symmetric polynomials.
Furthermore, $\mathcal{Q}$ can be identified by the set of the unramified points of the quotient variety $\mathbb{F}_q^m/{\rm{Sym}}_m$; see Section \ref{av}. The associated generalized reduced Datta-Johnsen code turns out to be equivalent to the evaluation code of the chosen linear (sub)system where the polynomials are evaluated on the set $\Delta$ of unramified points of the quotient variety $\mathbb{F}_q^m/{\rm{Sym}}_m$. In Section \ref{av}, an embedding of $\Delta$ in $AG(m,q)$ is described allowing us to use this model in our study.
%We also illustrate the difficulties in extending the results to the cases $m=2, t\ge 3$ and $m\ge 3,t\ge 2$.

A motivation of our investigation is to understand better how the fundamental parameters of the generalized Datta-Johnsen codes, especially their weight distributions, are related with enumerative questions concerning intersections of relevant objects in Finite geometry and Algebraic geometry over finite fields. We will see that such enumerative questions have already had satisfactory answers for $m=2,t=2$ and $q$ odd, and this allows us to work out this special case.  However, as long as $m\ge 3$ and $t\ge 2$, the study of the intersections of the involved algebraic varieties appears to be rather complicated, and is yet to carry out.

Our main results are described in Section \ref{ap}. In particular, a constructive proof of the following theorem is given.
\begin{theorem}
\label{main271224} For odd $q\ge 7$ there exist reduced generalized Datta-Johnsen codes  $\left[\frac{1}{2} q(q-1),3,D\right]$ whose minimum distance $D$ is at least $\frac{1}{2}\left(q^2-2q-2\sqrt{q}-7\right)$. The weights of the non-zero codewords fall into the interval $\left[\frac{1}{2} (q-1)-(\sqrt{q}+5),\frac{1}{2} (q-1)+(\sqrt{q}+3)\right]$.
\end{theorem} 
Computation for small values of $q=7,9$ shows that carefully chosen generalized Datta-Johnsen codes $\left[\frac{1}{2}q(q-1),3,d\right]$ have minimum distance $d$ equal to the optimal value minus 1.

Further investigations on the generalized Datta-Johnsen codes might  discover some more interesting features, especially about their locally recoverability; see \cite{Mi}.

\section{Background}
\subsection{Symmetric multivariable polynomials}
A polynomial $F\left(X_1,\ldots,X_m\right)\in \mathbb{K}\left[X_1,\ldots,X_m\right]$ in the indeterminates $X_1,\ldots,X_m$ with coefficients  over a field $\mathbb{K}$ is \emph{symmetric} if its invariant for any permutation $\pi$ on the index set  $\{1,\ldots,m\}$, i.e. $F\left(X_1,\ldots,X_m\right)=F\left(X_{\pi(1)},\ldots,X_{\pi(m)}\right)$.
For brevity, write $X=\left(X_1,\ldots,X_m\right)$. For $i=1,\ldots,m$, the $i$-th \emph{elementary symmetric polynomial} is
\[\sigma_{m}^i(X) = \sum_{1 \le j_1 < \ldots < j_i \le m} X_{j_1} \cdots X_{j_i}.\] In particular,
$\sigma_m^1(X)=X_1+\ldots+X_m$ and $\sigma_m^m(X)=X_1\cdots X_m$.

Let $f(Y)$ be any polynomial of degree $m$ in the unique indeterminate $Y$ with coefficients in a field $\mathbb{K}$. Let $y_1,\ldots,y_m$ be the (not necessarily distinct) roots of $f(Y)$ in an algebraic closure of $\mathbb{K}$. For brevity, write $y=(y_1,\ldots,y_m)$. Then the classical Viet\'e formula states
\begin{equation}\label{monic}
    f(Y)=\sum_{i=1}^m (-1)^i\sigma_m^i(y)Y^{m-i}+Y^m.
\end{equation}

For a polynomial $F(X_1,\ldots,X_m)\in \mathbb{K}[X_1,\ldots,X_m]$, substituting $X_i$ with the $i$-th symmetric polynomial provides a polynomial in $\mathbb{K}\left[X_1,\ldots,X_m\right]$:
$$G\left(X_1,\ldots,X_m\right)=F\left(\sigma_m^1(X),\ldots,\sigma_m^m(X)\right)$$ which is symmetric. The fundamental theorem on symmetric polynomials states that every symmetric polynomial $G\in \mathbb{K}\left[X_1,\ldots,X_m\right]$ arises in this way from a unique (not necessarily symmetric) polynomial $F\in \mathbb{K}\left[X_1,\ldots,X_m\right]$. This gives rise to a vector space monomorphism $\Phi_m$ from $\mathbb{K}[X_1,\ldots,X_m]$ onto its subspace $\mathbb{K}[X_1,\ldots,X_m]^s$ consisting of all symmetric polynomials. In particular, for every $\mathbb{K}$-subspace $\Sigma$ of $\mathbb{K}[X_1,\ldots,X_m]$ (also called linear system over $\mathbb{K}$), $\Phi_m(\Sigma)$ is a $\mathbb{K}$-subspace of symmetric polynomials, and the converse also holds.
\subsection{Affine varieties and curves over finite fields}
\label{av}
The affine space $AG(m,q)$ with a fixed coordinate system consists of all points $x=(x_1,\ldots,x_m)$ with $x_i\in \mathbb{F}_q$, and it is also an affine variety $\mathbb{F}_q^m$ defined over $\mathbb{F}_q$. The symmetric group ${\rm{Sym}}_m$ is a subgroup of $AGL(m,q)$ acting on the coordinates, i.e. if $g\in{\rm{Sym}}_m$ then $g(x)=(x_{i_1},\ldots,x_{i_m})$ where $[i_1,\ldots,i_m]$ is the permutation of $g$ on the coordinates.
This gives rise to a map $\psi_m$ from $\mathbb{F}_q^m$ to its quotient variety  $\mathbb{F}_q^m/{\rm{Sym}}_m$ whose points can be viewed as the orbits under the natural action of ${\rm{Sym}}_m$ on $AG(m,q)$. Let $\varphi_m$ be the rational map taking $x$ to the point $\varphi_m(x)=(\sigma_m^1(x),\ldots,\sigma_m^m(x))$. Then $\varphi_m$ commutes with every permutation on the coordinates, and hence $\varphi_m$ commutes with ${\rm{Sym}}_m$. For a point $y=\varphi_m(x)$, the fiber $\varphi_m^{-1}(y)$ consists of $x$ together with $\pi(x)$ where $\pi$ runs over the non-trivial elements of ${\rm{Sym}}_m$. In fact, if $\varphi_m(x)=\varphi_m(x')$ then $x_1,\ldots,x_m$ and $x_1',\ldots,x_m'$ are the roots of the same monic polynomial and Viet\'e's formula (\ref{monic}) yields that $x'$ and $x$ only differ by a permutation on the coordinates. Therefore, ${\rm{Im}}(\varphi_m)$ are the points of the quotient variety $\mathbb{F}_q^m/{\rm{Sym}}_m$, that is $\psi_m=\varphi_m$. Furthermore, the set $\mathcal{Q}$ defined in Introduction can be identified with the points of $\mathbb{F}_q^m/{\rm{Sym}}_m$.
Let $x$ be a point of $\mathbb{F}_q^m$. Then $x$ is a ramification point of $\varphi_m$, i.e. it contains two equal coordinates if and only if the monic polynomial $F(X)$ has a multiply root, and hence the discriminant of equation $F=0$ vanishes. Therefore, the set of all non-distinguished points $x$ is mapped by $\varphi_m$ into the discriminant variety $\Delta_m$.
\subsection{Upper bound on the number of points of a plane curve of \texorpdfstring{$PG(2,q)$}{Lg}}
Let $\mathcal{D}_u$ be a (possibly reducible) plane curve of $PG(2,q)$ of degree $u$ viewed as a curve of $PG(2,\bar{\mathbb{F}}_q)$ where $\bar{\mathbb{F}}_q$ stands for an algebraic closure of $\mathbb{F}_q$. 
Denote by $N_q(\mathcal{D}_u)$ the number of the points of $\mathcal{D}_u$ in $PG(2,q)$. It should be noticed that if $\mathcal{D}_u$ has some singular points then $N_q(\mathcal{D}_u)$ may not count the number of $\mathbb{F}_q$-rational points of a non-singular model of $\mathcal{D}_u$. In particular, the Hasse-Weil bound does not apply to $N_q(\mathcal{D}_u)$. The following upper bound is due to B. Segre \cite[Section 6, Theorems  I and II]{BS}:
\begin{equation}
\label{eq261224X}
|N_q(\mathcal{D}_u)|\le
\begin{cases}
(u-1)q+\left[\frac{1}{2} u\right]+1, &{\mbox{if $\mathcal{D}_u$ does not split into lines in $PG(2,\mathbb{F}_q$);}}\\
uq+1, &{\mbox{if $\mathcal{D}_u$ splits into lines over $\mathbb{F}_q$.}}\\
\end{cases}
\end{equation}
As a corollary,
\begin{equation}
\label{eq261224E}
{\mbox{$|N_q(\mathcal{D}_u)|\le uq+1$ for any plane curve $\mathcal{D}_u$ of degree $u$.}}
\end{equation}
In fact, if $\mathcal{C}_v$ is the product of all linear components of $\mathcal{C}_u$ defined over $\mathbb{F}_q$, and $\mathcal{C}_w$ the product of all the remaining components of $\mathcal{C}_u$, then  $|N_q(\mathcal{C}_u)|\le |N_q(\mathcal{C}_v)|+|N_q(\mathcal{C}_w)|$ and $u=v+w$. Hence the claim follows from (\ref{eq261224X}).

\subsection{Mutual positions of conics in planes of odd order}
\label{al} Let $\mathcal{C}$ be an irreducible conic in the projective plane $PG(2,q)$ of odd order $q$. The points of $PG(2,q)$ that are not on $\mathcal{C}$ are either external, which means that they lie on two tangents to $\mathcal{C}$, or internal, lying on no tangent to $\mathcal{C}$. Let $E_{\mathcal{C}}$ denote the set of all external points to $\mathcal{C}$. Then $E_{\mathcal{C}}$ has size $\frac{1}{2} q(q+1)$. The points of a second irreducible conic $\mathcal{D}$ in $PG(2,q)$  with respect to $\mathcal{C}$ fall into three subsets as the points of $\mathcal{D}$ other than those of $\mathcal{C}\cap \mathcal{D}$ are either internal to $\mathcal{C}$, or external to $\mathcal{C}$. The latter subset is $E_{\mathcal{C}}\cap \mathcal{D}$, and the possible sizes of $E_{\mathcal{C}}\cap \mathcal{D}$ are found in \cite{afkl}. There exist irreducible conics $\mathcal{D}$ whose points not on  $\mathcal{C}$ are all external, or all internal to $\mathcal{C}$. For such irreducible conics $\mathcal{D}$, assume that some point of $\mathcal{D}$ is external to $\mathcal{C}$. Then
%\begin{equation}
%\label{eq261224}
%|E_\mathcal{C}\cap \mathcal{D}|=
%\begin{cases} {\mbox{$q+1$, if $\mathcal{C}\cap\mathcal{D}=\emptyset$ and there is a cyclic projective group of order $q+1$ preserving  $\mathcal{C}$ and $\mathcal{D}$;}}\\
%{\mbox{$q$, if $|\mathcal{C}\cap\mathcal{D}|=1$ and there is a  projective group of order $q$  preserving $\mathcal{C}$ and $\mathcal{D}$;}}\\
%{\mbox{$q-1$, if $|\mathcal{C}\cap \mathcal{D}|=2$ and there is a cyclic projective group of order $q-1$ preserving  $\mathcal{C}$ and $\mathcal{D}$.}}
%\end{cases}
%\end{equation}
\begin{equation}
\label{eq261224}
|E_\mathcal{C}\cap \mathcal{D}|=
\begin{cases} q+1, &{\mbox{ if $\mathcal{C}\cap\mathcal{D}=\emptyset$ and there is a cyclic projective group of order}} \\ 
&{\mbox{ $q+1$ preserving  $\mathcal{C}$ and $\mathcal{D}$;}}\\
q, &{\mbox{ if $|\mathcal{C}\cap\mathcal{D}|=1$ and there is a  projective group of order $q$ }} \\
&{\mbox{ preserving $\mathcal{C}$ and $\mathcal{D}$;}}\\
q-1, &{\mbox{ if $|\mathcal{C}\cap \mathcal{D}|=2$ and there is a cyclic projective group of order }} \\
&{\mbox{ $q-1$ preserving  $\mathcal{C}$ and $\mathcal{D}$.}}
\end{cases}
\end{equation}
Clearly, $|E_{\mathcal{C}}\cap \mathcal{D}|=0$ if no point of $\mathcal{D}$ is external to $\mathcal{C}$. We count the total number of irreducible conics $\mathcal{D}$ in (\ref{eq261224}).  The projective group $PGL(2,q)$ preserving $\mathcal{C}$ has as many as $\frac{1}{2} q(q-1)$ cyclic subgroups of order $q+1$,
each has exactly $\frac{1}{2}(q-1)$ orbits which are conics contained in $E_{\mathcal{C}}$. Moreover, $PGL(2,q)$ has as many as $q+1$ subgroups of order $q$, and each has exactly $\frac{1}{2}(q-3)$ orbits which are conics contained in $E_{\mathcal{C}}$ minus a point of $\mathcal{C}$. Finally, $PGL(2,q)$ has as many as $\frac{1}{2} q(q+1)$ cyclic subgroups of order $q-1$ each has $\frac{1}{2}(q-3)$ orbits which are conics contained in $E_{\mathcal{C}}$ minus two points of $\mathcal{C}$. Therefore, the total number equals
\begin{equation}
\label{eq080125} \frac{1}{2} q(q-1) \frac{1}{2} (q-1)+(q+1)\frac{1}{2}(q-3)+\frac{1}{2} q(q+1) \frac{1}{2} (q-3)=\frac{1}{2} (q^3-q^2-3q-3).
\end{equation}
For any other irreducible conic $\mathcal{D}$, it is shown in \cite{afkl} that the size of $E_{\mathcal{C}}\cap \mathcal{D}$ can be computed from the number $N_q$ of $\mathbb{F}_q$-rational points of an elliptic curve defined over $\mathbb{F}_q$. The Hasse-Weil bound, see \cite[Theorem 9.18]{HKT}, $N_q-(q+1)\le 2\sqrt{q}$, yields the bound
\begin{equation}
\label{eq251224}
\frac{1}{2} (q-1)-\left(\sqrt{q}+3\right)\le |E_{\mathcal{C}}\cap \mathcal{D}|\le \frac{1}{2} (q-1)+\left(\sqrt{q}+3\right).
\end{equation}
For $q\le 19$, the exact value of the maximum size for $E_\mathcal{C}\cap \mathcal{D}$ is known, see \cite{afkl}. For conics $\mathcal{D}$ not considered in (\ref{eq261224}) these values are reported below.
 \begin{table}[!h]
\caption{Size of the largest intersection for cases not in (\ref{eq261224}) }
\label{smarc}
\hspace*{2mm}
\begin{tabular}{crrrrrrrrrrrrrr}
\hline
&&&&&&&&\\
&&&&&&&&\\[-6mm]
$q$      & 3 & 5 & 7 & 9 & 11 & 13  & 17 & 19 \\
$|E_\mathcal{C}\cap \mathcal{D}|$ & 1 & 3 & 5 & 7 & 9 & 10 & 13 & 13\\[2mm]
\hline
\end{tabular}
\end{table}

If $\mathcal{D}$ is reducible then it splits into two lines, say $\ell_1$ and $\ell_2$.
If $\ell_1=\ell_2$ then
\begin{equation}
\label{eqA251224}
|E_\mathcal{C}\cap \mathcal{D}|=
\begin{cases}
{\mbox{ $q$, if $\ell_1$ is tangent};} \\
{\mbox{ $\frac{1}{2}(q-1)$, if $\ell_1$ is secant};} \\
{\mbox{ $\frac{1}{2}(q+1)$, if $\ell_1$ is external}.} \\
\end{cases}
\end{equation}
If $\ell_1$ and $\ell_2$ are two distinct lines then
\begin{equation}
\label{eqB251224}
|E_\mathcal{C}\cap \mathcal{D}|=
\begin{cases}
2q-1, &{\mbox{if $\ell_1$ and $\ell_2$ are two distinct tangents};} \\
\frac{3}{2}(q-1)+1, &{\mbox{if $\ell_1$ is tangent and the common point of $\ell_1$ with $\ell_2$ is }}\\
& {\mbox{not on $\mathcal{C}$ };}\\
\frac{3}{2}(q-1), &{\mbox{if $\ell_1$ is tangent and $\ell_2$ is a secant such that their common }}\\
&{\mbox{point is on $\mathcal{C}$};}\\
q-1, &{\mbox{if both $\ell_1$ and $\ell_2$ are secants through either  internal point, }}\\
&{\mbox{or a point on $\mathcal{C}$};}\\
q-2, &{\mbox{if both $\ell_1$ and $\ell_2$ are secants through an external point};}\\
q, &{\mbox{if both $\ell_1$ and $\ell_2$ are external through an external point};} \\
q+1, &{\mbox{if both $\ell_1$ and $\ell_2$ are external through an internal point};}\\
q, &{\mbox{if  $\ell_1$ is secant and $\ell_2$ is external through an internal point};}\\
q-1, &{\mbox{if  $\ell_1$ is secant and $\ell_2$ is external through an internal point}.}
\end{cases}
 \end{equation}
 
\subsection{Linear codes}
A \emph{linear code} $ C $ of length $ n $ over a finite field $\mathbb{F}_q$ is a subspace of the vector space $\mathbb{F}_q^n$ over $\mathbb{F}_q$. The vectors in $C$ are the \emph{codewords}, and  if $C$ has dimension $k$ then it is a linear code of \emph{dimension} $k$. Fix a basis of $\mathbb{F}_q^n$. The \emph{weight} of a codeword is the number of its non-zero coordinates (entries). The \textit{Hamming distance} of two codewords $x,y\in C$ is the weight of $ x - y $. The \textit{minimum distance} $d$ of a code $ C $ is the minimum of distances of all two distinct codewords of $ C$ or, equivalently, the minimum weight of the non-zero vectors of $C$. A $[n,k,d]_q $\emph{-code} is a linear code with above parameters $n,k,d$.

One may ask for applications whether a given code is a ``good" one compared to others. Useful comparisons may be done with two further parameters, namely the \emph{relative distance} $\delta=d/n$ and the \emph{information} or \emph{dimension rate}  $R=k/n$. Codes with higher rates are considered to be better than codes with lower rates.
\subsubsection{Generalized Reed-Muller codes arising from multivariable polynomials over \texorpdfstring{$\mathbb{F}_q$}{Lg}}
The \emph{evaluation vector} of a polynomial $f\in \mathbb{F}_q[X_1,\ldots,X_m]$ on $(x_1,\ldots,x_m)$ with $x_i\in \mathbb{F}_q$ is the $m-tuple$ $(f(x_1),\ldots f(x_m))$, and it can be viewed as a vector in $\mathbb{F}_q^m$.
Let $x=(x_1,\ldots,x_m)$ with $x_i\in \mathbb{F}_q$. The $t$-th order Generalized Reed-Muller code is
\[GR_q(m,t):=\{(f(x) : x \in \mathbb{F}_q^m) \mid  f\in \mathbb{F}_q[x_1,\dots,x_m], \deg(f)\leq t\}\]
and it is a $\left[q^m, \binom{m+t}{m}, (1-\frac{t}{q})q^m\right]_q$ code.
\subsubsection{The Datta-Johnsen code}
\label{dj}The elementary symmetric polynomials $\sigma_m^i$ together with their $\mathbb{F}_q$-linear combinations form an ($m+1$)-dimensional $\mathbb{F}_q$-subspace in $\mathbb{F}_q[X_1,\ldots,X_m]$. Evaluating these polynomials on the set of all distinguished points in $\mathbb{F}_q^m$ (i.e. on all points with pairwise distinct coordinates in $\mathbb{F}_q^m$) is the Datta-Johnsen code $C_m$ introduced in \cite{datta2023codes}. For $m < q$, the code $C_m$ is a non-degenerate $[n,k,d]$ code, where $n =P(q,m)$ with
\begin{equation}
P(q,m)=
\begin{cases}
\binom{q}{m} m! \quad \text{if} \ \ m \le q,\\
\,\,\,\,0 \quad \text{otherwise},
\end{cases}
\end{equation} $k=m+1$ and $d =(q-m)P(q-1,m-1)$; see \cite[Proposition 3.2]{datta2023codes}. In \cite[Remark 3.3]{datta2023codes}, the authors pointed out that the distinguished points are partitioned into $\binom{q}{m}$ subsets each of which is an orbit of the symmetric group of degree $m$. Therefore, a smaller evaluation code $C_m'$ can be obtained by evaluating symmetric polynomials on an ordered set $\mathcal{Q}$ of representatives of those orbits.  For $m < q$ , $C_m'$ is a non-degenerate $[N, K, D]$ linear code where $N = \binom{q}{m}$, $K=m+1$ and $D=\binom{q}{m}- \binom{q-1}{m-1}$; see \cite[Proposition 3.4]{datta2023codes}.

\section{Evaluation codes of symmetric functions on distinguished points of the affine plane}
\label{ap}
 In this section, $m=2$. Therefore, $x=(x_1,x_2)\in\mathbb{F}_q^2$ and  $y_1=\sigma_2^1(x)=x_1+x_2$, $y_2=\sigma_2^2(x)=x_1x_2$. With this notation, $x_1=x_2$ if and only if $y_1=2x_1$ and $y_2=x_1^2$. Thus the discriminant variety $\Delta_2$ is the parabola $\mathcal{C}$ of equation $y_1^2-4y_2=0$. Take a distinguished point $P=(\xi_1,\xi_2)$. Then $\xi_1\neq \xi_2$, and  $\varphi_2(P)=(\sigma_2^1(P),\sigma_2^2(P))=(\xi_1+\xi_2,\xi_1\xi_2)$. Write $\eta_1=\sigma_2^1(P)$ and $\eta_2=\sigma_2^2(P)$. Now, look at the location of the point $Q=(\eta_1,\eta_2)$ with respect to the parabola $\mathcal{C}$. Clearly $Q\notin \mathcal{C}$. More precisely, we show that $\eta_1^2-4\eta_2$ is a non-zero square in $\mathbb{F}_q$.
 The polynomial in (\ref{monic}) for $m=2$ is $$F(X)=X^2-\eta_1X+\eta_2.$$ Since the roots $\xi_1$ and $\xi_2$ of the polynomial $F(X)$ are different and both defined over $\mathbb{F}_q$, its discriminant $\eta_1^2-4\eta_2$ is a non-zero square in $\mathbb{F}_q$. Conversely, if the point $Q=(\eta_1,\eta_2)$ is such that $\eta_1^2-4\eta_2$ is a non-zero square in $\mathbb{F}_q$ then the polynomial $F(X)=X^2-\eta_1X+\eta_2$ has two different roots $\xi_1$ and $\xi_2$ both in $\mathbb{F}_q$. The points $(\xi_1,\xi_2)$ and $(\xi_2,\xi_1)$ are two distinguished points in $\mathbb{F}_q^2$ and have the same $\varphi_2$-image. Moreover, the total number of distinguished points amounts $q(q-1)$. Therefore, that map $\varphi_2$ is a $2$-to-$1$ correspondence from the set of distinguished points onto a set of size $\frac{1}{2} q(q-1)$. Actually, this correspondence $\varphi_2$ becomes $1$-to $1$ it it is considered from a set $\mathcal{Q}$ of representatives, as defined in Section \ref{dj}. In geometric terms, $Im(\varphi_2)$ consists of all external points to the parabola  $\mathcal{C}$ in the affine plane $AG(2,\mathbb{F}_q)$.
 \subsection{Case of linear polynomials}
 \label{lp1}
Let $\Sigma_1$ be the linear system over $\mathbb{K}$ consisting of all linear polynomials in $x_1,x_2$. Then $\Phi_2(\Sigma_1)$ consists of all $\mathbb{F}_q$-linear combinations of $y_1$ and $y_2$. Evaluation of a symmetric polynomial $g\in \Phi_2(\Sigma_1)$ on the set of distinguished points in $\mathbb{F}_q^2$ can be carried out by evaluating the corresponding linear polynomial $f$ with $g=\Phi_m(f)$ on the set $E_\mathcal{C}$ of all affine external points to $\mathcal{C}$. In particular, the weight of the codeword of the Datta-Johnsen code $C_2'$ represented by the affine line $\ell$ of equation $f=0$ is equal to $\frac{1}{2} q(q-1)-(\ell\cap E_\mathcal{C})$. Here, $\ell\cap E_\mathcal{C}$ is either $q-1$, or $\frac{1}{2}(q-1)$, or $\frac{1}{2}(q-3)$ according as $\ell$ is a tangent to $\mathcal{C}$, or an external line to $\mathcal{C}$, or a secant to $\mathcal{C}$. Therefore, the weight distribution is $\{\frac{1}{2}(q-1)(q-2), \frac{1}{2}(q-1)^2, \frac{1}{2}(q^2-2q+3)\}$. Thus the reduced Datta-Johnsen code $C_2'$ is a $\left[\frac{1}{2} q(q-1),3,D\right]$ code which has minimum distance equal to $D=\frac{1}{2}(q-1)(q-2)$. Therefore, the Datta-Johnsen code has minimum distance $(q-1)(q-2)$ in accordance with \cite[Proposition 3.2]{datta2023codes} for the case $m=2$.
\subsection{Case of quadratic polynomials}
\label{lp2}
Let $\Sigma_2$ be the six-dimensional linear system over $\mathbb{K}$ consisting of all polynomials in $x_1,x_2$ of degree $\le 2$. After replacing $\Sigma_1$ by $\Sigma_2$ we may argue as in Section \ref{lp1}. The arising code has size $n=q(q-1)$ and dimension $6$ whose weight distribution depends on the  possible intersections between $E_\mathcal{C}$ and a conic $\mathcal{D}$ in the affine plane $AG(2,\mathbb{F}_q)$. These possibilities can be determined relying on Section \ref{al}. In particular, the maximum weight of the reduced generalized Datta-Johnsen code $C_2(6)'$\, is $2q-3$ and hence it is a $\left[\frac{1}{2} q(q-1),6,D\right]$ whose minimum distance $D$ equals $\frac{1}{2}(q-2)(q-3)$. Thus, the generalized Datta-Johnsen code $C_2(6)$ has minimum distance equal to $(q-2)(q-3)$. The linear system $\Sigma_2$ has linear subsystems $\Sigma_2(r)$ for any degree $r$ for $1\le r \le 5$. Each of them gives rise to a generalized Datta-Johnsen code of size $q(q-1)$ and dimension $r$ whose weight distribution and minimum distance can be estimated using the formulas (\ref{eq261224}), (\ref{eq251224}),  (\ref{eqA251224}) and (\ref{eqB251224}), and for smaller $q\le 19$ from Table \ref{smarc}. Here we limit ourselves to show two types of codes for $r=3$ with large minimum distances.

In the first type, the linear system $\Sigma_2(3)$ will contain no reducible conic whose components are lines defined over $\mathbb{F}_q$. To obtain such a linear system, look at the affine plane $AG(2,q^3)$ and its projective closure $PG(2,q^3)$ with homogeneous coordinates $(x:y:z)$. The projective  group $G=PGL(3,q)$ of $PG(2,q)$ can be viewed as a subgroup of $PGL(3,q^3)$. The action of $G$ on $PG(2,q^3)$ produces three point-orbits, namely $PG(2,q)$, the set of all points covered by lines of $PG(2,q)$ and the set $\Lambda$ of the remaining
points. Here
\begin{equation}
\label{eqA080125}
|\Lambda|=q^6+q^3+1-(q^2+q+1)-(q^2+q+1)(q^3-q)=q^6-q^5-q^4+q^3.
\end{equation}
Take a point $P=(a:b:c)\in\Lambda$  with its Frobenius images $P_1=(a^q:b^q:c^q)$ and $P_2=(a^{q^2}:b^{q^2}:c^{q^2})$. These points are the vertices of the triangle $PP_1P_2$ whose sides $\ell_1=PP_1$, $\ell_2=P_1P_2$ and $\ell_3=P_2P$ are disjoint from $PG(2,q)$. If $\ell_i$ has equation $\ell_i(x,y,z)=a_ix+b_iy+c_iz$ for $i=1,2,3$ then
$a_2=a_1^q,a_3=a_2^q,b_2=b_1^q,b_3=b_2^q,c_2=c_1^q,c_3=c_2^q$. This shows that the Frobenius image of $\ell_i$ is $\ell_{i+1}$ where the indices are taken$\pmod{3}$.

Define $\Sigma_2(3)$ to be the net (linear system of projective dimension $2$) of $PG(2,q^3)$ consisting of all conics through the points $P,P_1,P_2$. Clearly, $\Sigma_2(3)$ is generated by the three reducible conics, those of equations $\ell_1(x,y,z)\ell_2(x,y,z)=0$,  $\ell_2(x,y,z)\ell_3(x,y,z)=0$ and $\ell_3(x,y,z)\ell_1(x,y,z)=0$, respectively. Actually, $\Sigma_2(3)$ contains further reducible conics of equations $\ell t=0$ where $\ell$ coincides with a side of the triangle, and $t$ is any line through the opposite vertex of $\ell$. The line $t$ may have at most one point of $PG(2,q)$. Let $\bar{\Sigma}_2(3)$ be the set of all conics $\mathcal{C}_\lambda$ in $\Sigma_2(3)$ of equation $$\lambda\ell_1(x,y,z)\ell_2(x,y,z)+\lambda^q\ell_2(x,y,z)\ell_3(x,y,z)+\lambda^{q^2}\ell_3(x,y,z)\ell_1(x,y,z)=0$$ with $\lambda\in \mathbb{F}_{q^3}\setminus \{0\}$. The Frobenius image of $\mathcal{C}_\lambda$ is the conic of equation
$$\lambda^q\ell_2(x,y,z)\ell_3(x,y,z)+\lambda^{q^2}\ell_3(x,y,z)\ell_1(x,y,z)+\lambda^{q^3}\ell_1(x,y,z)\ell_2(x,y,z)=0.$$ Since $\lambda^{q^3}=\lambda$, this shows that $\bar{\Sigma}_2(3)$
consists of conics defined over $\mathbb{F}_q$, i.e. conics of $PG(2,q)$. The total number of conics in $\bar{\Sigma}_2(3)$ equals $(q^3-1)/(q-1)=q^2+q+1$, i.e. the number of points of $PG(2,q)$. It turns out that $\bar{\Sigma}_2(3)$ is a linear system of $PG(2,q)$ of projective dimension $2$. Moreover, $\bar{\Sigma}_2(3)$ contains no reducible conic. In fact, as we have already observed, any reducible conic in  $\Sigma_2(3)$ contains exactly one side of the triangle $PP_1P_2$ and hence it is not defined over $\mathbb{F}_q$ as the Frobenius map does not preserve any side of that triangle. We show that it is possible to choose the point $P$ such that none of the conics in (\ref{eq261224}) passes through $P$. From (\ref{eq080125}), the total number of such conics is smaller than $\frac{1}{2} q^3$. Each conic defined over $PG(2,q)$ has exactly $q^3-q$ points in $PG(2,q^3)$.
Therefore, the conics in (\ref{eq261224}) cover less than $\frac{1}{2} q^6$ points of $PG(2,q^3)$ outside $PG(2,q)$.
Since $\frac{1}{2} q^6<q^6-q^5-q^4+q^3$, they cannot cover all the points of $\Lambda$. Thus, any uncovered point in $\Lambda$ can be chosen for $P$.
From (\ref{eq251224}), the largest weight of the reduced generalized Datta-Johnson code $C_2(3)'$  does not exceed
$\frac{1}{2}(q+1)+\sqrt{q}+3$. Therefore, it is a $\left[\frac{1}{2} q(q-1),3,D\right]$ whose minimum distance $D$ is at least $\frac{1}{2} q(q-1)-\left(\frac{1}{2}(q+1)+\sqrt{q}+3\right)=\frac{1}{2}\left(q^2-2q-2\sqrt{q}-7\right)$. Thus, the generalized Datta-Johnsen code $C_2(3)$ has minimum distance not smaller than $q^2-2q-2\sqrt{q}-7$.

In the second type, the linear system $\Sigma_2(3)$ is defined over $\mathbb{F}_{q^2}$. In $PG(2,q^2)$ take a point $P$ not in $PG(2,q)$ together with its Frobenius image $P_1$. Then the line $\ell$ joining them is defined over $\mathbb{F}_q$, i.e. $\ell$ is a line of $PG(2,q)$.
Furthermore, take an internal point $P_2\in PG(2,q)$ to $\mathcal{C}$. Define $\Sigma_2(3)$ to be the linear system of conics in $PG(2,q^2)$ which passes through $P,P_1$ and $P_2$. Arguing as in the first example, it can be shown that the conics in $\Sigma_2(3)$ which are defined over $\mathbb{F}_q$ form a linear system $\bar{\Sigma}_2(3)$ in $PG(2,q)$. Actually, a direct presentation of $\bar{\Sigma}_2(3)$ is also possible. For this purpose, fix a non-square element $s$ of $\mathbb{F}_q$ and choose an element $i\in \mathbb{F}_{q^2}$ such that $i^2=s$. Then $i^q=-i$. Let $P=(i:1:0)$, $P_1=(-i:1:0)$ and $P_2=(0:-s:1)$. Then $P_2$ is an internal point to $\mathcal{C}$. With this notation,  $\bar{\Sigma}_2(3)$ consists of all conics $\mathcal{D}_{\alpha,\beta,\gamma}$ of equations
$$\alpha(x^2-sy^2)+2\beta xz+2\gamma yz+s(\alpha s^2+2\gamma )z^2=0.$$
%For $\mathcal{C}$, take the conic of equation $x^2-sy^2-z^2=0$, i.e. the parabola of equation $x^2-sy^2-1=0$ in the affine plane $AG(2,q)$ whose projective closure is $PG(2,q)$.  Then $P_2$ is an internal point to $\mathcal{C}$.
%With this notation,  $\bar{\Sigma}_2(3)$ consists of all conics $\mathcal{D}_{\alpha,\beta,\gamma}$ of equations
%$$\alpha(x^2-sy^2)+2\beta xz+2\gamma yz=0.$$
A straightforward computation shows that the determinant associated with $\mathcal{D}_{\alpha,\beta,\gamma}$ is equal to $-\alpha (\alpha s^2+\beta i+ \gamma)$ $(\alpha s^2-\beta i+ \gamma)$. Therefore, $\mathcal{D}_{\alpha,\beta,\gamma}$ is reducible if and only if either $\alpha=0$, or $\alpha s^2\pm \beta i+\gamma=0$. In the latter case, the reducible components are two (conjugate) lines defined over $\mathbb{F}_{q^2}$. For $\alpha=0$, one of the components is the line $\ell_\infty$ of equation $z=0$ and the other component is a line $t$ in $PG(2,q)$ passing through $P_2$. Therefore, if $\mathcal{D}$ is reducible then $E_\mathcal{C}\cap \mathcal{D}$ in $AG(2,q)$ only consists of the external points to $\mathcal{C}$ lying on $t$, and their number is at most $\frac{1}{2}(q-1)$. Moreover, none of the cases in (\ref{eq251224}) occurs as $\mathcal{D}_{\alpha,\beta,\gamma}$ contains an internal point to $\mathcal{C}$, namely $P_2$. We may conclude our discussion as in the first example. From (\ref{eq251224}), the largest weight of the reduced generalized Datta-Johnson code $C_2(3)'$  does not exceed
$\frac{1}{2}(q+1)+\sqrt{q}+3$. Therefore, it is a $\left[\frac{1}{2} q(q-1),3,D\right]$ whose minimum distance $D$ is at least $\frac{1}{2} q(q-1)-(\frac{1}{2}(q+1)+\sqrt{q}+3)=\frac{1}{2}(q^2-2q-2\sqrt{q}-7)$. Again, the generalized Datta-Johnsen code $C_2(3)$ has minimum distance not smaller than $q^2-2q-2\sqrt{q}-7$.
%\subsubsection{Case of cubic polynomials} The possibilities for the size of the intersection of $\mathcal{C}(E)$ with a plane cubic curve are not known.
\subsection{Case of polynomials of degree \texorpdfstring{$u\ge 3$}{Lg}} Unfortunately, an adaption for $3\le u <\frac{1}{2} q$ of the above idea, as developed in Section \ref{lp2}, works only partially for lack of the analog results on curves $\mathcal{D}$ of degree $u\ge 3$, quoted in Section \ref{al}. Nevertheless, some weaker general results can be obtained for $u\ge 3$ using (\ref{eq261224E}).
Let $\Sigma_u$ be the linear system of all plane curves $\mathcal{D}_u$ of degree $\leq u$. Then   the maximum weight of the reduced generalized Datta-Johnsen code $C_u'$\, is upper bounded $uq+1$ and hence it is a $\left[\frac{1}{2} q(q-1),\frac{1}{2} (u+1)(u+2),D\right]$ whose minimum distance $D$ is at least $\frac{1}{2} q(q-1)-uq-1$. Thus, the generalized Datta-Johnsen code $C_u$ has minimum distance at least $q(q-1-2u)-2$.

In the smallest case $u=3$, an analog of the first example for quadratic polynomials gives better results. With the previous notation, define $\bar{\Sigma}_3(4)$ to be the linear system of projective dimension $3$ of $PG(2,q^3)$ consisting of all cubics  $\mathcal{C}_\lambda$ of equation
\begin{equation}
    \label{eq18012025A}
\lambda\ell_1\ell_2^2+\lambda^q\ell_2\ell_3^2+\lambda^{q^2}\ell_3
\ell_1^2+\delta \ell_1 \ell_2\ell_3=0
\end{equation} where $\ell_i=\ell_i(x,y,z)$, $\lambda\in \mathbb{F}_{q^3}, \delta\in \mathbb{F}_q$, and $(\lambda,\delta)\ne (0,0)$. All these cubics are defined over $\mathbb{F}_q$ and their total number equals $(q^4-1)/(q-1)=q^3+q^2+q+1$, i.e. the number of points of $PG(3,q)$. Moreover,
the intersection multiplicity at the point $P_i$ of the triangle $PP_1P_2$ of $\mathcal{C}_\lambda\in \bar{\Sigma}_3(4)$ with a side $\ell_{i+1}$ is at least $2$, where $i=1,2,3, P_3=P$ and the indices are taken$\pmod{3}$. In fact, from the properties of the intersection multiplicity between plane curves, see \cite[Chapter 3]{HKT},  we have for $\lambda\ne 0$
\begin{equation}
\label{eq18012025}
I(P,\mathcal{C}_\lambda\cap \ell_1)=I(P,\ell_2\ell_3^2\cap \ell_1)=I(P,\ell_3^2\cap \ell_1)=2.
\end{equation}
We show that for $\bar{\Sigma}_3(4)$ contains no reducible cubic. Assume first that a cubic $\mathcal{C}$ in $\bar{\Sigma}_3(4)$ splits into an irreducible conic $C_2$ and a line $\ell$. Since $\ell$ is defined over $\mathbb{F}_q$, none of vertices of the triangle $PP_1P_2$ is incident with $\ell$. Therefore, $C_2$ passes through each vertex. This together with (\ref{eq18012025}) imply $I(P,C_2\cap \ell_1)=2$. Since $P_1\in C_2\cap \ell_1$, the Bézout theorem yields that $\ell_1$ is a component of $C_2$, a contradiction. Thus, $\mathcal{C}_\lambda$ splits into three lines. If one of these lines, say $\ell$ is defined over $\mathbb{F}_q$ then the choice of $P\in\Lambda$ rules out the possibility that a vertex of the triangle is incident with $\ell$. Then one of the other two lines passes through two vertices of the triangle, i.e. it coincides with one of the sides $\ell_i$. Since $\mathcal{C}_\lambda$ is defined over $\mathbb{F}_q$, the Frobenius images of $\ell$ are also components of $\mathcal{C}_\lambda$. Thus $\lambda=0$.

Notice that the projective dimension $3$ of $\bar{\Sigma}_3(4)$ is best possible. In fact, by \cite[Theorem 1.3]{AGR24}, if $\Sigma$ is a linear system of cubics with projective dimension at least $4$, then $\Sigma$ has an $\mathbb{F}_q$-member which is reducible over $\mathbb{F}_q$.

Now, let $\mathcal{C}_\lambda\in \bar{\Sigma}_3(4)$ of affine equation obtained from (\ref{eq18012025A}) by putting $z=1$. Also, write $\ell_i(x,y)=\ell_i(x,y,1)$. Moreover replace $y$ by $\frac{1}{4}(x^2-t^2)$. The arising affine polynomial has equation
\begin{multline}
f(x,t)=\lambda\ell_1(x,\frac{1}{4}(x^2-t^2))\ell_2(x,\frac{1}{4}(x^2-t^2))^2+\\
\lambda^q\ell_2(x,\frac{1}{4}(x^2-t^2))\ell_3(x,\frac{1}{4}(x^2-t^2))^2+
\lambda^{q^2}\ell_3(x,\frac{1}{4}(x^2-t^2))\ell_1(x,\frac{1}{4}(x^2-t^2))^2+\\
\delta \lambda_1(x, \frac{1}{4}(x^2-t^2))\lambda_2(x,\frac{1}{4}(x^2-t^2))\lambda_3(x,\frac{1}{4}(x^2-t^2)).
\end{multline}
Then a point in $Q(\xi,\eta)$ of $AG(2,q)$ is a common point of $E_\mathcal{C}$ and $\mathcal{C}_\lambda$ if and only if $f(\xi,\tau)=0$ where $\eta=\frac{1}{4}(\xi^2-\tau^2)$,  $\xi,\tau\in \mathbb{F}_q$  and $\tau\ne 0$; in other words the point $Q'(\xi,\tau)$ with $\tau \ne 0$ is an affine point of the plane curve $\mathcal{F}$ of equation $f(x,t)=0$. Therefore, to compute the weight-distribution of the reduced generalized Datta-Johnsen code $C_3(4)'$ it is necessary to know the number of points of $\mathcal{F}$ in $AG(2,q)$. Unfortunately, the exact value of the number $N_q(\mathcal{U})$ of points that a given plane algebraic curve $\mathcal{U}$ of degree $6$ defined over $\mathbb{F}_q$ possesses, is largely unknown. To estimate the minimum distance of $C_3(4)'$, upper bounds on $N_q(\mathcal{U})$ are required. If either $\mathcal{U}$ is absolutely irreducible, or it has an absolutely irreducible component defined over $\mathbb{F}_q$, the Hasse-Weil bound, the Serre bound, or the St\"ohr-Voloch bound may be useful, especially for larger values of $q$. Useful criteria for low dimensional linear systems to consist of absolutely irreducible plane curves are found in \cite{AGR24}.   On the other hand, for smaller values of $q$, an exhaustive computation can be carried out. Table \ref{smarc1} shows the computational results for $q=5,7,9,11$ and it makes a comparison with the conics cases and with the known bounds on the minimum distances for codes with same length and dimension.

\begin{table}[!h]
\caption{}
\label{smarc1}
\begin{tabular}{lcllcc}
            & \textbf{Dim} & \textbf{Curves} & \textbf{Max distances} & \textbf{Upper bound} &  \\ \hline\hline
$q = 5$     & 3                  & Type 1 conics   & 5, 6, 7                & 7                    & from \cite{Grassl}      \\
Length = 10 & 3                  & Type 2 conics   & 5                      & 7                    & from \cite{Grassl}      \\
            & 4                  & Cubics          & 3--6             & 6                    & from \cite{Grassl}      \\ \hline
$q = 7$     & 3                  & Type 1 conics   & 14, 15, 16             & 17                   & from \cite{Grassl}      \\
Length = 21 & 3                  & Type 2 conics   & 15                     & 17                   & from \cite{Grassl}      \\
            & 4                  & Cubics          & 11--14                  & 16                   & Delsarte LP      \\ \hline
$q = 9$     & 3                  & Type 1 conics   & 27--30                  & 31                   & from \cite{Grassl}      \\
Length = 36 & 3                  & Type 2 conics   & 28                     & 31                   & from \cite{Grassl}      \\
            & 4                  & Cubics          & 23--27                  & 30                   & from \cite{Grassl}      \\ \hline
$q = 11$    & 3                  & Type 1 conics   & 44--47                  & 49                   & Delsarte LP      \\
Length = 55 & 3                  & Type 2 conics   & 45                     & 49                   & Delsarte LP      \\
            & 4                  & Cubics          & 40--44                  & 48                   & Delsarte LP      \\
\textbf{}   &                    &                 &                        &                      &
\end{tabular}
\end{table}

For $q=7,9,11$, the extremal cases can be obtained as follows.
Fix a primitive element $w$ of $\mathbb{F}_{q^3}$, and let $w_1=w^{q-1}$. For $1\le k\le q^2+q$, define the line $\ell_1$ to be of affine equation $w_1^kx+w_1y+1=0$. Let $\mathcal{C}^{(k)}_\lambda$ be the irreducible plane cubic of equation (\ref{eq18012025A}), and denote by $N(k)$ the maximum of $|E_\mathcal{C}\cap \mathcal{C}^{(k)}_\lambda|$ where $\lambda$ ranging on all non-zero elements of $\mathbb{F}_{q^3}$. Then, the minimum distance of the reduced generalized Datta-Johnsen code $C_3(4)'$ is equal to $\frac{1}{2} q(q-1)-\frac{1}{2} N(k)$; see Table \ref{smarc2} for the computed values when $q=5,7,9,11$.

 \begin{table}[!ht]
\caption{}
\label{smarc2}
\hspace*{2mm}
\begin{tabular}{crrrrrrrr}
\hline
&&&&&&&&\\
&&&&&&&&\\[-6mm]
$q$      & 5 & 7 & 9 & 11  \\ [1mm]
$k$      & 11 & 18 & 23 & 17 \\ [1mm]
$N(k)$   & 10 & 14 & 18 & 22 \\  [1mm]
$\frac{1}{2} q(q-1)-\frac{1}{2} N(k)$\,\,\,\, & 5 & 14 & 27 & 44 \\[2mm]
\hline
\end{tabular}
\end{table}

\section{Evaluation codes of symmetric functions on distinguished points of the affine space}
In this section, $m=3$. For $x=(x_1,x_2,x_3)\in\mathbb{F}_q^3$, let $y_1=\sigma_3^1(x)=x_1+x_2+x_3,y_2=\sigma_3^2=x_1x_2+x_2x_3+x_1x_3$ and $y_3=\sigma_3^3=x_1x_2x_3$. %Let compute the discriminant variety $\Delta_3$.
For a point $P=(\xi_1,\xi_2,\xi_3)$ we have $\varphi_3(P)=(\sigma_3^1(P),\sigma_3^2(P),\sigma_3^3(P))=(\xi_1+\xi_2+\xi_3,\xi_1\xi_2+\xi_2\xi_3+\xi_1\xi_3,\xi_1\xi_2\xi_3)$. Write $\eta_i=\sigma_3^i(P)$ for $i=1,2,3$. The polynomial  (\ref{monic}) reads
\begin{equation}
    \label{eq30112024}X^3-\eta_1X^2+\eta_2X-\eta_3.
\end{equation}
%Its roots are $\xi_1,\xi_2,\xi_3$.
If $p>3$, the classical Cardano's formula for the roots holds true in characteristic $p$. Therefore, if $p>3$ and $\xi$ is a root of the polynomial (\ref{eq30112024}) then
$$ \xi=\{b + [b^2 + (c-a^2)^3]^{1/2}\}^{1/3} + \{b -[b^2 + (c-a^2)^3]^{1/2}\}^{1/3}+a$$
where
$$a =\frac{1}{3}\eta_1,\, b = \frac{1}{27}\eta_1^3 - \frac{1}{6}(\eta_1\eta_2-3\eta_3),\,   c = \frac{1}{3}\eta_3.$$

Therefore, $\xi\in \mathbb{F}_q$ if and only if $\xi^q-\xi=0$, that is,
\begin{equation}
\label{eqA30112024}
\begin{array}{lll}
\Big(\{\big(\frac{1}{27}\eta_1^3-\frac{1}{6}(\eta_1\eta_2-3\eta_3)
+[\frac{1}{27}\eta_1^3-\frac{1}{6}(\eta_1\eta_2-3\eta_3)]^2-
(\frac{1}{3}\eta_3-\frac{1}{9}\eta_1^2)^3\big)^{1/2}\}^{1/3}+\\
\{\big(\frac{1}{27}\eta_1^3-\frac{1}{6}(\eta_1\eta_2-3\eta_3)
-[\frac{1}{27}\eta_1^3-\frac{1}{6}(\eta_1\eta_2-3\eta_3)]^2-
(\frac{1}{3}\eta_3-\frac{1}{9}\eta_1^2)^3\big)^{1/2}\}^{1/3}\Big)^q-\\
\{\big(\frac{1}{27}\eta_1^3-\frac{1}{6}(\eta_1\eta_2-3\eta_3)
+[\frac{1}{27}\eta_1^3-\frac{1}{6}(\eta_1\eta_2-3\eta_3)]^2-
(\frac{1}{3}\eta_3-\frac{1}{9}\eta_1^2)^3\big)^{1/2}\}^{1/3}+\\
\{\big(\frac{1}{27}\eta_1^3-\frac{1}{6}(\eta_1\eta_2-3\eta_3)
-[\frac{1}{27}\eta_1^3-\frac{1}{6}(\eta_1\eta_2-3\eta_3)]^2-
(\frac{1}{3}\eta_3-\frac{1}{9}\eta_1^2)^3\big)^{1/2}\}^{1/3}=0.
\end{array}
\end{equation}
Thus, $\varphi_3$ maps $AG(3,\mathbb{F}_q)$ to the set $\Delta$ of all points $P=(\eta_1,\eta_2,\eta_3)$ whose coordinates satisfy Equation (\ref{eqA30112024}). In other words, the points $(\eta_1,\eta_2,\eta_3)$ of the quotient variety $\mathbb{F}_{q^3}/{\rm{Sym}}_3$ are those satisfying Equation (\ref{eqA30112024}).
The polynomial (\ref{monic}) has a multiple root, i.e. $P$ is not a distinguished point, if and only if the system
\begin{equation*}
\begin{cases}
    X^3-\eta_1X^2+\eta_2X-\eta_3=0, \\
    3X^2-2\eta_1X+\eta_2=0
    \end{cases}
\end{equation*}
has a solution.
From the Sylvester determinant
\begin{equation*}
\begin{vmatrix}
 1& -\eta_1 & \eta_2 & -\eta_3 & 0\\
    0 & 1& -\eta_1 & \eta_2 & -\eta_3\\
    3 & -2\eta_1 & \eta_2 & 0 & 0 \\
    0 & 3 & -2\eta_1 & \eta_2 & 0 \\
    0 & 0 & 3 & -2\eta_1 & \eta_2
\end{vmatrix}
=-\eta_1^2\eta_2^2+4\eta_2^3+4\eta_1^3\eta_3+27\eta_2^2-18\eta_1\eta_2\eta_3,
\end{equation*}
the discriminant variety $\Delta_3$ is the (affine) surface $\mathcal{F}$ of equation
$y_1^2y_2^2-4y_2^3-4y_1^3y_3-27y_2^2+18y_1y_2y_3=0$
which shows that the set of all non-distinguished points $x$ is mapped into $\Delta_3$. Therefore, $\varphi_3$ maps the set of all distinguished points onto $\Delta\setminus \Delta_3$. Fix an order of the distinguished points $\mathcal{Q}=\{Q_1,\ldots,Q_n\}$ with $n=\frac{1}{6}q(q-1)(q-2)$.
Evaluating the linear polynomials in $\mathbb{F}_q[X_1,X_2,X_3]$ on $\Delta\setminus \Delta_3$ gives the reduced Datta-Johnsen code $C_3'$. Its weight distributions of $C_3'$ determine (and determined by) the sizes of the possible plane sections of the quotient variety $\mathbb{F}_{q^3}/{\rm{Sym}}_3$ and the determinental variety. If quadratic polynomials in $\mathbb{F}_q[X_1,X_2,X_3]$
considered, then an analog connection between weight distribution and sizes of sections of $\mathbb{F}_{q^3}/{\rm{Sym}}_3$ and the determinental variety occurs where the sections are cut out by quadrics. Looking back to (\ref{eqA30112024}) the study of such sections appears to be out of reach.

\section{Acknowledgments} 
B. Gatti, G. Korchm\'aros and G. Schulte  have been partially supported by the Italian National Group for Algebraic and Geometric Structures and their Applications (GNSAGA - INdAM). V. Pallozzi Lavorante has been  
partially supported by NSF grant number 2127742

%%===========================================================================================%%
%% If you are submitting to one of the Nature Portfolio journals, using the eJP submission   %%
%% system, please include the references within the manuscript file itself. You may do this  %%
%% by copying the reference list from your .bbl file, paste it into the main manuscript .tex %%
%% file, and delete the associated \verb+\bibliography+ commands.                            %%
%%===========================================================================================%%

%\bibliography{GattiKorchmarosOthers}% common bib file
%% if required, the content of .bbl file can be included here once bbl is generated
%%\input {GattiKorchmarosOthers.bbl}

\end{document}